\documentclass[12pt,reqno]{amsproc}%
\usepackage{color,anysize,float,CJK}
\usepackage{citeref}
\usepackage{txfonts,esint}
\DeclareSymbolFont{CM}{OMX}{cmex}{m}{n}
\DeclareMathSymbol{\sumop}{\mathop}{CM}{"50}
\renewcommand{\sum}{\sumop}
\usepackage[colorlinks,urlcolor=blue,citecolor=blue,linkcolor=blue]{hyperref}
\usepackage{amsmath}
\usepackage{amsfonts}
\usepackage{amssymb}
\usepackage{graphicx}%

\marginsize{25mm}{25mm}{15mm}{15mm}
\renewcommand{\bibitempages}[1]{}

\newtheorem{thm}{Theorem}[section]

\newtheorem{lem}[thm]{Lemma}
\newtheorem{prop}[thm]{Proposition}
\theoremstyle{definition}

\theoremstyle{remark}
\newtheorem{rem}[thm]{\it Remark}

\numberwithin{equation}{section}
\allowdisplaybreaks

\begin{document}
\title[Standing waves for 6-superlinear Chern-Simons-Schr\"{o}dinger systems]{Standing waves for 6-superlinear Chern-Simons-Schr\"{o}dinger systems with indefinite potentials}

\author{Shuai Jiang$^*$\and Shibo Liu}
\dedicatory{School of Mathematical Sciences, Xiamen University,
	Xiamen 361005, China}

\thanks {${^\star}$This work was supported by NSFC (Grant No. 12071387 and 11971436).\\
\rule{4ex}{0ex}2020 Mathematics Subject Classification. 35J20, 35J50, 35J10.\\
\rule{3ex}{0ex}$^*$Corresponding author.\\
\rule{4ex}{0ex}\emph{E-mail address:} jiangshuai0915@163.com (S. Jiang), liusb@xmu.edu.cn (S. Liu).}
\maketitle
\vspace{-5mm}

\begin{abstract}
	In this paper we consider 6-superlinear Chern-Simons-Schr\"{o}dinger systems. In contrast to most studies, we consider the case where the potential $V$ is indefinite so that the Schr\"{o}dinger operator $-\Delta +V$ possesses a finite-dimensional negative space. We obtain nontrivial solutions for the problem via Morse theory.
\vskip1ex

\end{abstract}

\section{Introduction}
In this paper, we consider the following Chern-Simons-Schr\"{o}dinger system (CCS system) in $H^1(\mathbb{R}^{2})$:
\begin{equation}
	\left\{
	\begin{array}[c]{ll}%
		-\Delta u+V(x)u+A_0u+\sum_{j=1}^{2}{A_j}^2u=f(x,u),\\
	\partial_1A_0=A_2\vert u\vert^2,\quad\partial_2A_0=-A_1\vert u\vert^2,\\
	\partial_1A_2-\partial_2A_1=-\frac{1}{2}u^2,\quad\partial_1A_1+\partial_2A_2=0.
	\end{array}
	\right.
	\label{eq:K}%
\end{equation}
 where $V\in C(\mathbb{R}^{2})$ is potential and $f\in C(\mathbb{R}^{2},\mathbb{R})$ is the nonlinearity. The   $\left(\mathbf{C}\mathbf{C}\mathbf{S}\right)$ system describes the nonrelativistic thermodynamic behavior of large number of particles in an electromagnetic field.

 The $\left(\mathbf{C}\mathbf{C}\mathbf{S}\right)$ system (\ref{eq:K}) arises when we are looking for standing waves for the following nonlinear Schr\"{o}dinger system
 \begin{equation}
 	\left\{
 \begin{array}[c]{ll}%
 	iD_0\phi+\left(D_1D_1+D_2D_2\right)\phi+f(x,\phi)=0,\\
 	\partial_0A_1-\partial_1A_0=-\text{Im}\left(\bar{\phi}D_2\phi\right),\\
 	\partial_0A_2-\partial_2A_0=\text{Im}\left(\bar{\phi}D_1\phi\right),\\
 	\partial_1A_2-\partial_2A_1=-\frac{1}{2}\vert\phi\vert^2,
 		\end{array}
 \right.\label{eq:12}
 \end{equation}
where $i$ denotes the imaginary unit, $\partial_0=\frac{\partial}{\partial_t}$, $\partial_1=\frac{\partial}{\partial_{x_1}}$, $\partial_2=\frac{\partial}{\partial_{x_2}}$, for $(t,x_1,x_2)\in \mathbb{R}^{1+2}$,  $\phi:\mathbb{R}^{1+2}\rightarrow\mathbb{C}$ is the complex scalar field, $A_\mu:\mathbb{R}^{1+2}\rightarrow\mathbb{R}$ is the gauge field. The associated covariant differential operators are given by
\[
D_\mu:=\partial_\mu+iA_\mu,\qquad\mu=0,1,2.
\]
System (\ref{eq:12}) proposed in \cite{MR1084552,MR1056846} consists of the Schr\"{o}dinger equation augmented by the gauge field $A_\mu$. This feature of the model is important for the study of the high-temperature superconductor, fractional quantum Hall effect and Aharovnov-Bohm scattering.

We suppose that the gauge field satisfies the Coulmb gauge condition $\partial_0A_0+\partial_1A_1+\partial_2A_2=0$, and $A_\mu(x,t)=A_\mu(x)$, $\mu=0,1,2.$ Then we deduce that $\partial_1A_1+\partial_2A_2=0$. Moreover,
standing waves for (\ref{eq:12}) are obtained through the ansatz $\phi=u(x)e^{i\omega t}$, $f(x,ue^{i\omega t})=f(x,u)e^{i\omega t}$, $\omega>0$, resulting in
\begin{equation}
	\left\{
\begin{array}[c]{ll}%
	-\Delta u+\omega u+A_0u+\sum_{j=1}^{2}{A_j}^2u=f(x,u),\\
	\partial_1A_0=A_2\vert u\vert^2,\quad\partial_2A_0=-A_1\vert u\vert^2,\\
	\partial_1A_2-\partial_2A_1=-\frac{1}{2}u^2,\quad \partial_1A_1+\partial_2A_2=0.
\end{array}
\right.
\label{eq:13}%
\end{equation}
Here the components $A_1$ and $A_2$ in system (\ref{eq:13}) can be represented by solving the elliptic equation
\[
\Delta A_1=\partial_2\left(\frac{\vert u\vert^2}{2}\right)\quad\text{and}\quad\Delta A_2=-\partial_2\left(\frac{\vert u\vert^2}{2}\right),
\]
which provide
\[
A_1=A_1[u](x)=\frac{x_2}{2\pi\vert x\vert^2}*\left(\frac{|u|^2}{2}\right)=\frac{1}{2\pi}\int_{\mathbb{R}^2}\frac{x_2-y_2}{\vert x-y\vert^2}\frac{\vert u(y)\vert^2}{2}\text{d}y,
\]
\[A_2=A_2[u](x)=\frac{x_1}{2\pi\vert x\vert^2}*\left(\frac{\vert u\vert^2}{2}\right)=-\frac{1}{2\pi}\int_{\mathbb{R}^2}\frac{x_1-y_1}{\vert x-y\vert^2}\frac{\vert u(y)\vert^2}{2}\text{d}y,
\]
where $*$ denotes the convolution.

Similarly,
\[
\Delta A_0=\partial_1(A_2\vert u\vert^2)-\partial_2(A_1\vert u\vert^2),
\]
which gives the following representation of the component $A_0$:
\[
A_0=A_0[u](x)=\frac{x_1}{2\pi\vert x\vert^2}*(A_2\vert u\vert^2)-\frac{x_2}{2\pi\vert x\vert^2}*(A_1\vert u\vert^2).
\]

The $\left(\mathbf{C}\mathbf{C}\mathbf{S}\right)$ system (\ref{eq:K}) has attracted considerable attention  in recent decades, which can be seen in \cite{MR2948224,MR3494398,MR3619082,MR3173176,MR4188316,MR4000150,MR3924547}  and the references therein. We emphasize that in all these papers, the authors only considered the case where the Schr\"{o}dinger operator $-\Delta+V$ is positive definite. In this case, the quadratic part of the variational functional $\Phi$ given in (\ref{eq:2.1}) is positively definite, the zero function $u=0$ is a local minimizer of $\Phi$ and the mountain pass theorem \cite{MR0370183} can be applied. However, when the potential $V$ is negative somewhere so that the quadratic part of $\Phi$ is indefinite, the zero function $u=0$ is no longer a local minimizer of $\Phi$, the mountain pass theorem is not applicable anymore. For stationary NLS equations
\[
-\Delta u+V(x)u=f(x,u)
\]
with indefinite Schr\"{o}dinger operator $-\Delta+V$, one usually applies the linking theorem to get solution, see e.g.\cite{MR1751952,MR2957647}. For system (\ref{eq:K}), it seems hard to verify the linking geometry due to the nonnegative terms involving $A_i^2u^2$ in \eqref{eq:2.1}, which prevent the functional $\Phi$ to be nonpositive on the negative space of the Schr\"{o}dinger operator.  Hence, the classical linking theorem \cite[Lemma 2.12]{MR1400007} is also not applicable.

For this reason, there are very few results about (\ref{eq:K}) with indefinite potential.  It seems that \cite{MR4271204} is the only work devoted to this situation. To overcome these difficulties, and the difficulty that the Sobolev embedding
\[
H^1(\mathbb{R}^2)\hookrightarrow L^2(\mathbb{R}^2)
\]
is not compact, it is assumed in  \cite{MR4271204} that, roughly speaking, $V$ is coercive so that the related Sobolev space is compactly embedded into     $ L^2(\mathbb{R}^2)$. Then the local linking theory of Li and Willem \cite{MR1312028} is applied to get critical point of $\Phi$.

Motivated by the above observation and \cite{MR3656292} on Schr\"{o}dinger-Poisson systems (see also \cite{MR3945610}), in this paper, we will consider the case that $V$ is bounded, so that the above-mentioned compact embedding may not be true. From now on all integrals are taken over $\mathbb{R}^2$ except stated explicitly.

Now we are ready to state our assumptions on $V$ and $f$.

\begin{itemize}
	\item[$(V)$] $V\in C(\mathbb{R}^{2})$ is a bounded function such that the quadratic form
	\[
	Q(u)=\frac{1}{2}\int\left(\vert \nabla u\vert
	^{2}+ V(x)u^{2}\right)
	\]
	is non-degenerate and the negative space of $Q$ is finite-dimensional.
	\item[$(f_1)$] $f\in C(\mathbb{R}^2\times\mathbb{R})$ satisfies
	\[
	 \lim_{\left\vert
		t\right\vert \rightarrow0}\frac{f(x,t)}{t}=0\text{,}\qquad\lim_{\left\vert t\right\vert \rightarrow\infty}\frac{f(x,t)}{e^{\mu t^2}}=0
	\]
	for any $\mu>0$, where $F(x,t)=\int_{0}^{t}f(x,\tau)\text{d}\tau$.
	\item[$(f_2)$] For $(x,t)\in \mathbb{R}^2\times\mathbb{R}\setminus\{0\}$ we have $0<6F(x,t)\leq tf(x,t)$, moreover, for almost all $x\in\mathbb{R}^2$
	\begin{equation}
		\lim_{\left\vert
			t\right\vert \rightarrow\infty}\frac{F(x,t)}{t^6}=+\infty.\label{eq:1.4}
	\end{equation}
	\item[$(f_3)$] One of the following conditions is satisfied:
	\begin{enumerate}
	\item[$(f_{31})$] there exist $C_0>0$ and $\nu\in(0,6)$ satisfying $
	F(x,t)\ge C_0\vert t\vert^\nu$ for all $t\in\mathbb{R}$;
	\item[$(f_{32})$] for some $\delta>0$, $
	F(x,t)\leq0$ for all $\vert t\vert\leq\delta$.
	\end{enumerate}
		
	\item[$(f_4)$] For any $r>0$, we have
		\[
		\lim_{\left\vert x\right\vert \rightarrow\infty}\underset{0<\left\vert t\right\vert\leq r}{\sup}\left\vert\frac{f(x,t)}{t}\right\vert=0.
		\]
		\item[$(f_5)$] For some $s>2$, $p,q>1$ we have $a\in L^\infty(\mathbb{R}^2) \cap L^{p}(\mathbb{R}^2), b\in L^\infty(\mathbb{R}^2) \cap L^{q}(\mathbb{R}^2)$ such that
	\begin{equation}
		\left\vert f(x,t)\right\vert\leq a(x)\vert t\vert+b(x)\vert t\vert^{s-1}\text{.}\label{eq:1.5}
	 \end{equation}
\end{itemize}
We emphasize that in $(f_5)$, the exponents $p$ and $q$ can be chosen arbitrarily from $(1,\infty)$, see Remark \ref{rr0}.

Now we are ready to state the main results of this paper.
\begin{thm}
\label{th:1}Suppose that $(V)$, $(f_1)$, $(f_2)$, $(f_3)$ and $(f_4)$ are satisfied. Then system \eqref{eq:K} has a nontrivial solution.
\end{thm}
\begin{thm}
	\label{th:2}Suppose that $(V)$, $(f_1)$, $(f_2)$, $(f_3)$ and $(f_5)$ are satisfied. Then system \eqref{eq:K} has a nontrivial solution.
\end{thm}
\begin{rem}
	\label{re:1}
	As we have mentioned before, neither the mountain pass lemma nor the linking theorem can be applied to our functional $\Phi$. It turns out that $\Phi$ has a local linking at the origin. Unfortunately, at present all critical point theorems involving local linking require the functional to satisfy global compactness condition. The role of our condition $(f_4)$ or $(f_5)$ is to ensure such compactness.
\end{rem}
The paper is organized as follows. In Section 2 we prove that the $(PS)$ sequences of $\Phi$ are bounded and $\Phi$ satisfies the $(PS)$ condition. In Section 3 we will prove Theorem \ref{th:1} by applying Morse theory. For this purpose after recalling some concepts and results in infinite-dimensional Morse theory \cite{MR1196690}, we will compute the critical group of $\Phi$ at infinity and then give the proof of Theorem \ref{th:1}. Finally, after investigating the compactness of the operator $K^{\prime}$ (see Lemma \ref{lf:1}), we use Morse theory again to prove Theorem \ref{th:2} in Section 4.
\section{Palais-Smale Condition}
Throughout this paper we always denote $X=H^1(\mathbb{R}^2)$. Under the assumptions $(V)$ and $(f_1)$, similar to Kang and Tang \cite{MR4271204} we can show that the functional $\Phi:X\rightarrow\mathbb{R}$,
\begin{equation}
\Phi(u)=\frac{1}{2}\int\left(\vert\nabla u\vert
^{2}+ V(x)u^{2}+A^2_1(u)u^2+A^2_2(u)u^2\right)
-\int F(x,u). \label{eq:2.1}%
\end{equation}
is well defined and of class $C^1$. The derivative of $\Phi$ is given by
\[
\left\langle \Phi^{\prime}(u),v\right\rangle =\int\nabla u\cdot\nabla v+\int V(x)uv+\int\left[\left(A^2_1(u)+A^2_2(u)\right)uv+A_0uv\right]-\int f(x,u)v.
\]
for $u,v\in X$. Consequently, critical points of $\Phi$ are weak solutions of system \eqref{eq:K}.

To study the functional $\Phi$, it will be convenient to rewrite the quadratic part $Q$ in a simpler form. It is well known that, if $(V)$ holds, then there exists an equivalent norm $\Vert\cdot\Vert$ on $X$ such that
\[
Q(u)=\frac{1}{2}\left(\Vert u^+\Vert^2-\Vert u^-\Vert^2\right),
\]
where $u^\pm$ is the orthogonal projection of $u$ on $X^\pm$ being $X^\pm$ the positive/negative space of $Q$. Using this new norm, $\Phi$ can be rewritten as
\begin{equation}
	\Phi(u)=\frac{1}{2}\left(\Vert u^+\Vert^2-\Vert u^-\Vert^2\right)+\frac{1}{2}\int\left(A^2_1(u)u^2+A^2_2(u)u^2\right)
	-\int F(x,u). \label{eq:2.2}
\end{equation}
By simple calculation (also see \cite{MR3619082,MR3924547}), we obtain, for any $u\in X$,
\begin{equation}
\left\langle \Phi^{\prime}(u),u\right\rangle =\Vert u^+\Vert^2-\Vert u^-\Vert^2+3\int\left(A^2_1(u)u^2+A^2_2(u)u^2\right)
-\int f(x,u)u. \label{eq:2.3}
\end{equation}

Next, we recall the following properties of the terms involving $A_0$, $A_1$, $A_2$.
\begin{lem}[\cite{MR4271204}]
	\label{l:2.1}  There is a constant $a_1>0$ such that for all $u\in X$ we have
	\[
	0\leq\int\left(A^2_1(u)+A^2_2(u)\right)u^2\leq a_1\left\Vert u\right\Vert^6.
	\]
\end{lem}
\begin{lem}[{\cite[Proposition 2.1]{MR3619082}}]
	\label{l:2.2} Let $1<r<2$ and $\frac{1}{r}-\frac{1}{t}=\frac{1}{2}$. If $u\in X$, then
\begin{align}
	&\vert A_0(u)\vert_t\leq C\vert u\vert^2_{2r}\vert u\vert^2_4\nonumber,\\
	&\vert A^2_i(u)\vert_t\leq C\vert u\vert^2_{2r}\nonumber,
\end{align}
	where $i=1,2$ and $\vert\cdot\vert_t$ is the $L^t$-norm.
\end{lem}

\begin{lem}
 \label{l:2.3}If $(V)$, $(f_1)$ and $(f_2)$ hold, then all $(PS)$ sequences of $\Phi$ are bounded.
\end{lem}
\begin{proof}
	Let $\{u_n\}$ be a $(PS)$ sequence of $\Phi$,  that is,
	\[
	\sup_n\vert\Phi(u_{n})\vert<\infty,\qquad\Phi^{\prime}\left(u_{n}\right)\rightarrow{0}.
	\]
	It suffices to show that $\{u_n\}$ is bounded.
	Suppose $\{u_n\}$ is unbounded, we may assume $\Vert u_n\Vert\rightarrow\infty.$ Let $v_n=\Vert u_n\Vert^{-1}u_n$. Then
	\[
	v_n=v^+_n+v^-_n\rightharpoonup v=v^++v^-\in X,\quad v^\pm_n,v^\pm\in X^\pm.
	\]

	If $v=0$, then $v^-_n\rightarrow v^-=0$ because dim $X^-<\infty$. Since
	\[
	\Vert v^+_n\Vert^2+\Vert v^-_n\Vert^2=1,
	\]
	for $n$ large enough we have
	\[
	\Vert v^+_n\Vert^2-\Vert v^-_n\Vert^2\ge\frac{1}{2}.
	\]
	Therefore, by assumption $(f_2)$, we deduce that for $n$ large enough,
	\begin{align}
		1+\sup_n\vert\Phi(u_{n})\vert+\Vert u_n\Vert
		&\ge\Phi(u_n)-\frac{1}{6}\left\langle\Phi^{\prime}(u_n),u_n\right\rangle\nonumber\\
		&=\frac{1}{3}\Vert u_n\Vert^2\left(	\Vert v^+_n\Vert^2-\Vert v^-_n\Vert^2\right)+\int\left(\frac{1}{6}f(x,u_n)u_n-F(x,u_n)\right)\nonumber\\
		&\ge\frac{1}{6}\Vert u_n\Vert^2\nonumber,
	\end{align}
contradicting $\Vert u_n\Vert\rightarrow\infty$.

If $v\neq0$. Then the set $\Theta=\{v\neq0\}$ has positive Lebesgue measure. For $x\in\Theta$ we have $\vert u_n(x)\vert\rightarrow\infty$ and
\begin{equation}
\frac{F(x,u_n(x))}{\Vert u_n\Vert^6}=\frac{F(x,u_n(x))}{u^6_n(x)}v^6_n(x)\rightarrow+\infty\label{eq:2.4},
\end{equation}
thanks to (\ref{eq:1.4}).
By Fatou lemma we deduce from (\ref{eq:2.4}) that
\begin{equation}
\int\frac{F(x,u_n)}{\Vert u_n\Vert^6}\ge\int_{v\neq0}\frac{F(x,u_n)}{\Vert u_n\Vert^6}\rightarrow+\infty.\label{eq:2.5}
\end{equation}
It follows from Lemma \ref{l:2.1} that
\begin{align}
\frac{1}{\Vert u_n\Vert^6}\int F(x,u_n)
&=\frac{\Vert u^+_n\Vert^2-\Vert u^-_n\Vert^2}{2\Vert u_n\Vert^6}+\frac{1}{2\Vert u_n\Vert^6}\int \left(A^2_1(u_n)u_n^2+A^2_2(u_n)u_n^2\right)-\frac{\Phi(u_n)}{\Vert u_n\Vert^6}\nonumber\\
&\leq\frac{a_1}{2}+1\nonumber,
\end{align}
a contradiction to (\ref{eq:2.5}). Therefore $\{u_n\}$ is bounded in $X$.
\end{proof}
To get a convergent subsequence of the $(PS)$ sequence, we need some compact properties of operators involving $A_j$, $j=$0,1,2. Firstly, we need to investigate the $C^1$-functional $\mathcal{N}:X\rightarrow \mathbb{R}$,
\[
\mathcal{N}(u)=\frac{1}{2}\int\left(A^2_1(u)u^2+A^2_2(u)u^2\right).
\]
It is known that the derivative of $\mathcal{N}$ is given by
\[
\left\langle \mathcal{N}^{\prime}(u),v\right\rangle=\int\left[\left(A^2_1(u)+A^2_2(u)\right)uv+A_0(u)uv\right],\quad u,v\in X.
\]
\begin{lem}
	\label{l:2.4} The functional $\mathcal{N}$ is weakly lower semi-continuous, its derivative $\mathcal{N}^{\prime}: X\rightarrow X^*$ is weakly sequentially continuous, where $X^*=H^{-1}(\mathbb{R}^2)$ is the dual space of $X=H^1(\mathbb{R}^2)$.
\end{lem}
\begin{proof}
	Let $\{u_n\}$ be a sequence in $X$ such that $u_n\rightharpoonup u$ in $X$, we need to show
	\[
	\mathcal{N}(u)\leq \varliminf \mathcal{N}(u_n),\qquad\left\langle \mathcal{N}^{\prime}(u_n),\phi\right\rangle\rightarrow\left\langle \mathcal{N}^{\prime}(u),\phi\right\rangle,
	\]
	for all $\phi\in X$.

	Since $u_n\rightharpoonup u$ in $X$, up to a subsequence, by the compactness of the embedding $X\hookrightarrow L^2_{\rm loc}(\mathbb{R}^2)$, we have
	\[
	u_n\rightarrow u\quad\text{in}~~L^2_{\rm loc}(\mathbb{R}^2),\qquad u_n(x)\rightarrow u(x)\quad \text{a.e.\ in}~~\mathbb{R}^2.
	\]
	According to Wan and Tan \cite[Proposition 2.2]{MR3619082}, for $j=1,2$ we have $A^2_j(u_n)\rightarrow A^2_j(u)$ a.e. in $\mathbb{R}^2$. Moreover,
	by the Fatou lemma,
	\[
	\mathcal{N}(u)=\frac{1}{2}\int\left(A^2_1(u)u^2+A^2_2(u)u^2\right)\leq\frac{1}{2}\varliminf\int\left(A^2_1(u_n)u_n^2+A^2_2(u_n)u_n^2\right)=\varliminf \mathcal{N}(u_n).
	\]
	Hence  $\mathcal{N}$ is weakly lower semi-continuous.
	
	To prove the weak continuity of $\mathcal{N}^{\prime}$, we observe that
	\begin{align}
	\left\langle\mathcal{N}^{\prime}(u_n)-\mathcal{N}^{\prime}(u),\phi\right\rangle
	=&\int\left(A^2_1(u_n)u_n\phi-A^2_1(u)u\phi\right)+\int\left(A^2_2(u_n)u_n\phi-A^2_2(u)u\phi\right)\nonumber\\
	&+\int\left(A_0(u_n)u_n\phi-A_0(u)u\phi\right).\label{eq:2.6}
\end{align}
Since
\begin{equation}
	\int\left(A^2_j(u_n)u_n\phi-A^2_j(u)u\phi\right)=\int A^2_j(u_n)(u_n-u)\phi+\int\left(A^2_j(u_n)-A^2_j(u)\right)u\phi,\quad j=1,2.\label{eq:2.7}
\end{equation}
By Lemma \ref{l:2.2}, H\"{o}lder inequality and continuous embedding yield that
	\begin{align}
	\int\vert A^2_j(u_n)(u_n-u)\vert^2
	&\leq\vert A^2_j(u_n)\vert^2_4\vert u_n-u\vert^2_4\nonumber\\
	&\leq C\vert u_n\vert^4_{\frac{8}{3}}\vert u_n-u\vert^2_4\nonumber\\
	&\leq C\Vert u_n\Vert^4\Vert u_n-u\Vert^2\leq C.\nonumber
\end{align}
Combining $u_n\rightarrow u$ a.e.\ in $\mathbb{R}^2$ and $A^2_j(u_n)\rightarrow A^2_j(u)$ a.e.\ in $\mathbb{R}^2$, we have $ A^2_j(u_n)(u_n-u)\rightharpoonup 0$ in $L^2(\mathbb{R}^2)$. Thus
 \begin{equation}
 \int A^2_j(u_n)(u_n-u)\phi\rightarrow0.\label{eq:2.8}
 \end{equation}
Similarly,
\begin{equation}
	\int\left(A^2_j(u_n)-A^2_j(u)\right)u\phi\rightarrow0.\label{eq:2.9}
\end{equation}
On the other hand
\begin{equation}
	\int\left(A_0(u_n)u_n\phi-A_0(u)u\phi\right)=\int A_0(u_n)(u_n-u)\phi+\int\left(A_0(u_n)-A_0(u)\right)u\phi\label{eq:2.10}.
\end{equation}
By Lemma \ref{l:2.2}, H\"{o}lder inequality and continuous embedding yield that
\begin{align}
	\int\vert A_0(u_n)(u_n-u)\vert^2
	&\leq\vert A_0(u_n)\vert^2_4\left\vert u_n-u\right\vert^2_4\nonumber\\
	&\leq C\vert u_n\vert^4_{\frac{8}{3}}\vert u_n\vert^4_4\vert u_n-u\vert^2_4\nonumber\\
	&\leq C\Vert u_n\Vert^8\Vert u_n-u\Vert^2\leq C.\nonumber
\end{align}
Combining $u_n\rightarrow u$ a.e.\ in $\mathbb{R}^2$, we have $ A_0(u_n)(u_n-u)\rightharpoonup0$ in $L^2(\mathbb{R}^2)$. Thus
\begin{equation}
	\int A_0(u_n)(u_n-u)\phi\rightarrow0.\label{eq:2.11}
\end{equation}
Similarly,
\begin{equation}
	\int\left(A_0(u_n)-A_0(u)\right)u\phi\rightarrow0.\label{eq:2.12}
\end{equation}
From (\ref{eq:2.6})-(\ref{eq:2.12}), for $\phi\in X$ we have
\[
\left\langle \mathcal{N}^{\prime}(u_n),\phi\right\rangle\rightarrow\left\langle \mathcal{N}^{\prime}(u),\phi\right\rangle.
\]
Therefore, we have proved that $\mathcal{N}^{\prime}$ is weakly sequencetially continuous.
\end{proof}
\begin{lem}
	\label{l:2.5}Let $u_n\rightharpoonup u$ in $X$. Then
	\[
	\varliminf\int\left[\left(A^2_1(u_n)+A^2_2(u_n)\right)u_n(u_n-u)+A_0(u_n)u_n(u_n-u)\right] \ge0.
	\]
\end{lem}
\begin{proof}
	Applying Lemma \ref{l:2.4}, we have
	\[
	\varliminf\mathcal{N}(u_n)\ge\mathcal{N}(u),\qquad \text{lim}\left\langle \mathcal{N}^{\prime}(u_n),u\right\rangle=\left\langle \mathcal{N}^{\prime}(u),u\right\rangle.
	\]
	Therefore,
	\begin{align}
		&\hspace{-1cm}\varliminf\int\left[\left(A^2_1(u_n)+A^2_2(u_n)\right)u_n(u_n-u)+A_0(u_n)u_n(u_n-u)\right]\nonumber\\
		=&\varliminf\int\left[3\left(A^2_1(u_n)+A^2_2(u_n)\right)u_n^2-\left(A^2_1(u_n)+A^2_2(u_n)\right)u_nu-A_0(u_n)u_nu\right]\nonumber\\
		=&\varliminf\left(6\mathcal{N}(u_n)-\left\langle \mathcal{N}^{\prime}(u_n),u\right\rangle\right)\nonumber\\
		\ge&6\mathcal{N}(u)-\left\langle \mathcal{N}^{\prime}(u),u\right\rangle=0.\qedhere\nonumber
	\end{align}
\end{proof}
\begin{lem}
\label{l:2.6}If $(V)$, $(f_1)$, $(f_2)$ and $(f_4)$ hold, then the functional
$\Phi$ satisfies the $(PS)$ condition, that is, any $(PS)$ sequence $\{u_{n}\}\subset
X$ possesses a convergent subsequence.
\end{lem}
\begin{proof}
	Let $\{u_n\}$ be a $(PS)$ sequence. We know from Lemma \ref{l:2.3} that $\{u_n\}$ is bounded in $X$. Up to a subsequence we may assume $u_n\rightharpoonup u$ in $X$. We have
	\[
	\int\left(\nabla u_n\cdot\nabla u+V(x)u_nu\right)\rightarrow\int\left(\vert \nabla u\vert
	^{2}+ V(x)u^{2}\right)=\Vert u^+\Vert^2-\Vert u^-\Vert^2.
	\]
	Consequently
	\begin{align}
	o(1)=&\left\langle \Phi^{\prime}(u_n),u_n-u\right\rangle\nonumber\\
	=& \int\left[\nabla u_n\cdot\nabla \left({u_n-u}\right)+V(x)u_n\left(u_n-u\right)\right]\nonumber\\
	&+\int\left[\left(A^2_1(u_n)+A^2_2(u_n)\right)u_n(u_n-u)+A_0(u_n)u_n(u_n-u)\right]-\int f(x,u_n)(u_n-u)\nonumber\\
	=&\int\left(|\nabla u_n|^2+V(x)u_n^2\right)-\int\left(\nabla u_n\cdot\nabla u+V(x)u_nu\right)\nonumber\\
	&+\int\left[\left(A^2_1(u_n)+A^2_2(u_n)\right)u_n(u_n-u)+A_0(u_n)u_n(u_n-u)\right]-\int f(x,u_n)(u_n-u)\nonumber\\
	=&\left(\Vert u_n^+\Vert^2-\Vert u_n^-\Vert^2\right)-\left(\Vert u^+\Vert^2-\Vert u^-\Vert^2\right)\nonumber\\
	&+\int\left[\left(A^2_1(u_n)+A^2_2(u_n)\right)u_n(u_n-u)+A_0(u_n)u_n(u_n-u)\right]-\int f(x,u_n)(u_n-u)+o(1).\nonumber
\end{align}
	  We have $u^-_n\rightarrow u^-$ and $\Vert u_n^-\Vert=\Vert u^-\Vert$ because dim$X^-<\infty$. Collecting all infinitesimal terms, we obtain
	\begin{align}
    \Vert u_n^+\Vert^2-\Vert u^+\Vert^2
	= o(1)&+\int f(x,u_n)(u_n-u)\nonumber\\
	&-\int\left[\left(A^2_1(u_n)+A^2_2(u_n)\right)u_n(u_n-u)+A_0(u_n)u_n(u_n-u)\right].\label{eq:2.13}	
\end{align}
Using the condition $(f_4)$, according to \cite[p.29]{MR2038142} we have
   \[
	\varlimsup\int f(x,u_n)(u_n-u)\leq0.
   \]
We deduce from Lemma \ref{l:2.5} and (\ref{eq:2.13}) that
\begin{align}
	&\hspace{-1cm}\varlimsup\left(\Vert u_n^+\Vert^2-\Vert u^+\Vert^2\right)\nonumber\\
	=&\varlimsup\left(\int f(x,u_n)(u_n-u)-\int\left[\left(A^2_1(u_n)+A^2_2(u_n)\right)u_n(u_n-u)+A_0(u_n)u_n(u_n-u)\right]\right)\nonumber\\
	=&\varlimsup\int f(x,u_n)(u_n-u)-\varliminf\int\left[\left(A^2_1(u_n)+A^2_2(u_n)\right)u_n(u_n-u)+A_0(u_n)u_n(u_n-u)\right]\nonumber\\
	\leq&\varlimsup\int f(x,u_n)(u_n-u)\leq0.\nonumber
\end{align}
Combining this with the weakly lower semi-continuity of the norm functional $u\mapsto \left\Vert u\right\Vert$, we obtain
\[
 \Vert u^+\Vert\leq\varliminf \Vert u^+_n\Vert\leq\varlimsup \Vert u^+_n\Vert\leq \Vert u^+\Vert.
\]
Therefore $ \Vert u^+_n\Vert\rightarrow \Vert u^+\Vert$. Remembering $ \Vert u^-_n\Vert\rightarrow \Vert u^-\Vert$, we get $ \Vert u_n\Vert\rightarrow \left\Vert u\right\Vert$. Thus $u_n\rightarrow u$ in $X$.
\end{proof}

\section{Critical groups and the proof of Theorem \ref{th:1}}

Having established the $(PS)$ condition for $\Phi$, we are now ready to present the proof of Theorem \ref{th:1}. We start by recalling some concepts and results from infinite-dimensional Morse theory (see e.g., Chang \cite{MR1196690} and Mawhin and Willem \cite[Chapter 8]{MR982267}).

Let $X$ be a Banach space, $\varphi: X\rightarrow \mathbb{R}$ be a $C^1$ functional, $u$ be an isolated critical point of $\varphi$ and $\varphi(u)=c$. Then
\[
C_q(\varphi,u):=H_q(\varphi_c,\varphi_c\setminus\{0\}),\qquad q\in\mathbb{N}=\{0,1,2,\ldots\},
\]
is called the $q$-th critical group of $\varphi$ at $u$, where $\varphi_c:=\varphi^{-1}(-\infty,c]$ and $H_*$ stands for the singular homology with coefficients in $\mathbb{Z}$.

If $\varphi$ satisfies the $(PS)$ condition and the critical values of $\varphi$ are bounded from below by $\alpha$, then following Bartsch and Li \cite{MR1420790}, we define the $q$-th critical group of $\varphi$ at infinity by
\[
C_q(\varphi,\infty):=H_q(X,\varphi_\alpha),\qquad q\in\mathbb{N}.
\]
Due to the deformation lemma, it is well known that the homology on the right hand side does not depend on the choice of $\alpha$.
\begin{prop}[{\cite[Proposition 3.6]{MR1420790}}]
	\label{prop:3.1} If $\varphi\in C^1(X,\mathbb{R})$ satisfies the $(PS)$ condition and $C_\ell(\varphi,0)\neq C_\ell(\varphi,\infty)$ for some $\ell\in\mathbb{N}$, then $\varphi$ has a nonzero critical point.
\end{prop}
\begin{prop}[{\cite[Theorem 2.1]{MR1110119}}]
	\label{prop:3.2}  Suppose $\varphi\in C^1(X,\mathbb{R})$ has a local linking at $0$ with respect to the decomposition $X=Y\oplus Z$, i.e., for some $\varepsilon>0$,
	\begin{align}
		&\varphi(u)\leq0 \quad\text{for}~~u\in Y\cap B_\varepsilon,\nonumber\\
		&\varphi(u)>0 \quad\text{for}~~u\in(Z\setminus\{0\})\cap B_\varepsilon,\nonumber
	\end{align}
where $B_\varepsilon=\{u\in X\mid\Vert u\Vert\leq\varepsilon\}$. If $\ell={\rm dim}~Y<\infty$, then $C_\ell(\varphi,0)\neq0$
\end{prop}
To investigate $C_*(\Phi,\infty)$, using the idea of \cite[Lemma 3.3]{MR3656292} we will prove the following lemma.
\begin{lem}
	\label{l:3.3}If $(V)$, $(f_1)$ and $(f_2)$ hold, there exists $A>0$ such that, if $\Phi(u)\leq-A$, then
	\[
	\frac{\rm d}{{\rm d}t}\Bigg|_{t=1}\Phi(tu)<0.
	\]
\end{lem}
\begin{proof}
	Otherwise, there exists a sequence $\{u_n\}\subset X$ such that $\Phi(u_n)\leq-n$ but
	\begin{equation}
	\left\langle \Phi^{\prime}(u_n),u_n\right\rangle=	\frac{\rm d}{{\rm d}t}\Bigg|_{t=1}\Phi(tu_n)\ge0.\label{eq:3.1}
\end{equation}
Consequently,
\begin{align}
	2\left(\Vert u_n^+\Vert^2-\Vert u_n^-\Vert^2\right)
	&\leq2\left(\Vert u_n^+\Vert^2-\Vert u_n^-\Vert^2\right)+\int\left[f(x,u_n)u_n-6F(x,u_n)\right]\nonumber\\
	&=6\Phi(u_n)-\left\langle \Phi^{\prime}(u_n),u_n\right\rangle\leq-6n\label{eq:3.2}
\end{align}
Let $v_n=\Vert u_n\Vert^{-1}u_n$ and $v^\pm_n$ be the orthogonal projection of $v_n$ on $X^\pm$. Then up to a subsequence $v^-_n\rightarrow v^-$ for some $v^-\in X^-$, because dim $X^-<\infty$.

If $v^-\neq 0$, then $v_n\rightharpoonup v$ in $X$ for some $v\in X\setminus\{0\}$. By assumption $(f_2)$ we have
\[
\frac{f(x,t)t}{t^6}\ge\frac{6F(x,t)}{t^6}\rightarrow+\infty,\qquad\text{as}~~t\rightarrow\infty.
\]
Thus, similar to the proof of $(\ref{eq:2.5})$, we obtain
\begin{equation}
	\frac{1}{\Vert u_n\Vert^6}\int f(x,u_n)u_n\rightarrow+\infty.\label{eq:3.3}
\end{equation}
Now, using Lemma \ref{l:2.1} we have a contradiction
\begin{align}
	0\leq\frac{\left\langle \Phi^{\prime}(u_n),u_n\right\rangle}{\Vert u_n\Vert^6}
	&=\frac{1}{\Vert u_n\Vert^6}\left(\left(\Vert u_n^+\Vert^2-\Vert u_n^-\Vert^2\right)+3\int\left(A^2_1(u_n)+A^2_2(u_n)\right)u^2_n-\int f(x,u)u\right)\nonumber\\
	&\leq 1+3a_1-\frac{1}{\Vert u_n\Vert^6}\int f(x,u_n)u_n\rightarrow-\infty.\nonumber
\end{align}
Hence, we must have $v^-=0$. But $\Vert v^+_n\Vert^2+\Vert v^-_n\Vert^2=1$, we deduce $\Vert v^+_n\Vert\rightarrow1$. Now for large $n$ we have
\[
\Vert u^+_n\Vert=\Vert u_n\Vert\Vert v^+_n\Vert\ge\Vert u_n\Vert\Vert v^-_n\Vert=\Vert u^-_n\Vert
\]
a contradiction to $(\ref{eq:3.2})$.
\end{proof}
\begin{lem}
	\label{l:3.4} $C_q(\Phi,\infty)=0$ for all $q\in\mathbb{N}$.
\end{lem}
\begin{proof}
		Let $B=\{v\in X\mid \Vert u\Vert\leq1\}$, $S=\partial{B}$ be the unit sphere in $X$, and $A>0$ be the number given in Lemma \ref{l:3.3}. Without loss of generality, we may assume that
		\[
		-A<\underset{0<\left\Vert u\right\Vert\leq2}{\inf}\Phi(u).
		\]
		Using $(\ref{eq:1.4})$, it is easy to see that for any $v\in S$
		\begin{align}
		\Phi(sv)&=\frac{s^2}{2}\left(\Vert v^+\Vert^2-\Vert v^-\Vert^2\right)+3s^2\int\left(A^2_1(sv)+A^2_2(sv)\right)v^2-\int F(x,sv)\nonumber\\
		&=s^6\left\{\frac{\Vert v^+\Vert^2-\Vert v^-\Vert^2}{2s^4}+\frac{3}{s^4}\int\left(A^2_1(sv)+A^2_2(sv)\right)v^2-\int\frac{F(x,sv)}{s^6}\right\}\rightarrow-\infty,\nonumber
	\end{align}
as $s\rightarrow+\infty$. Therefore,  for $v\in S$ there is $s_v>0$ such that $\Phi(s_vv)=-A$.

Using Lemma \ref{l:3.3} and the implicit function theorem, as in the proof of \cite[Lemma 3.4]{MR3656292} it can be shown that such $s_v$ is uniquely determined by $v$ and $T:v\mapsto s_v$ is continuous on $S$. Using the continuous function $T$ it is standard (see \cite{MR1094651}) to construct a deformation from $X\setminus B$ to the level set $\Phi_{-A}=\Phi^{-1}(-\infty,-A\rbrack$, and deduce
\[
C_q(\Phi,\infty)=H_q(X,\Phi_{-A})\cong H_q(X,X\setminus B)=0,\qquad\text{for all}~~q\in\mathbb{N}.\qedhere
\]
\end{proof}

\begin{proof}
	[\indent Proof of Theorem \ref{th:1}]
	From the assumption $(f_1)$, there exists $\varepsilon>0$ and $C_\varepsilon>0$ such that, for every $\mu>0$,
	\[
	\vert F(x,u)\vert\leq\varepsilon\vert u\vert^2+C_\varepsilon\vert u\vert^\xi(e^{\mu u^2}-1)\quad \text{for all}~~\xi>2.
	\]
	Using the conditions $(V)$, $(f_1)$ and $(f_3)$, similar to \cite[Lemmas 3.4 and 3.5]{MR4271204}, it is easy to see  that
    \[
			\Phi(u)=\frac{1}{2}\left(\Vert u^+\Vert^2-\Vert u^-\Vert^2\right)+o(\Vert u\Vert^2)\quad\text{as $\Vert u\Vert\rightarrow0$.}
	\]
	Hence, there exists $\varepsilon>0$ such that $\Phi$ is positive on $(X^+\setminus\{0\})\cap B_\varepsilon$, and negative on $(X^-\setminus\{0\})\cap B_\varepsilon$. That is, $\Phi$ has a local linking with respect to the decomposition $X=X^-\oplus X^+$. Therefore Proposition \ref{prop:3.2} yields
	\[
	 C_\ell(\varphi,0)\neq0,
	\]
	where $\ell=\text{dim}~X^-$. By Lemma \ref{l:3.4}, $C_\ell(\varphi,\infty)=0$. Applying Proposition \ref{prop:3.1}, we see that $\Phi$ has a nonzero critical point. The proof of Theorem \ref{th:1} is completed.
\end{proof}
	\section{Proof of Theorem \ref{th:2}}
	To prove Theorem \ref{th:2}, we need to recover the $(PS)$ condition with the help of condition $(f_5)$. Therefore, we need the following lemma.
	\begin{lem}
		\label{lf:1} Suppose that $f:\mathbb{R}^2\times\mathbb{R}\rightarrow\mathbb{R}$ is continuous and satisfies $(f_5)$.		
		 For the functional $K:X\rightarrow \mathbb{R},$
		\[
		K(u)=\int F(x,u)
		\]
		is well defined and of class $C^1$ with
		\[
		\left\langle K^{\prime}(u),\phi\right\rangle=\int f(x,u)\phi,\quad\forall\phi\in X.
		\]
		Moreover, $K^{\prime}$ is compact.
	\end{lem}
\begin{proof}
	From \eqref{eq:1.5} we have
	\[
	|f(x,t)|\le |a|_\infty|t|+|b|_\infty|t|^{s-1}\text{,}
	\]
	so it is well knwon that $K$ is well defined and of class $C^1$.
	
To show that $K':X\to X^*$ is compact, let $u_n\rightharpoonup u$ in $X$.	
Because $a\in L^{p}(\mathbb{R}^2), b\in  L^{q}(\mathbb{R}^2)$, for any $\varepsilon>0$, there exists $R>0$ such that
	  \begin{equation}
	  \int_{B^{\,\rm c}_R}\vert a\vert^p<\varepsilon^p,\quad\int_{B^{\,\rm c}_R}\vert b\vert^q
	  <\varepsilon^q,\label{eq:43}
	  \end{equation}
where $B_R$ is the ball in $\mathbb{R}^2$ with radius $R>0$ centering at the origin and $B^{\,\rm c}_R=\mathbb{R}^2\setminus B_R$.

For $\phi\in X$, $\|\phi\|=1$, using (\ref{eq:1.5}), (\ref{eq:43}) and the H\"{o}lder inequality and noting that $u,\phi\in L^\gamma(B^{\,\rm c}_R)$ for any $\gamma>2$, we have
	 \begin{align}
	  \int_{B^{\,\rm c}_R}\vert f(x,u)\vert\vert\phi\vert
	  &\leq\int_{B^{\,\rm c}_R}\vert a\vert\vert u\vert\vert\phi\vert+\int_{B^{\,\rm c}_R}\vert b\vert\vert u\vert^{s-1}\vert\phi\vert\nonumber\\
	  &\leq[ a]_{p}[ u]_{2p/(p-1)}[\phi]_{2p/(p-1)}+[ b]_q[ u]^{s-1}_{2q(s-1)/(q-1)}[\phi]_{2q/(q-1)}\nonumber\\
	  &<\varepsilon\left([ u]_{2p/(p-1)}[\phi]_{2p/(p-1)}+[ u]^{s-1}_{2q(s-1)/(q-1)}[\phi]_{2q/(q-1)}\right)\nonumber\\
	  &\le M\varepsilon\text{,}
     \end{align}
where $[\cdot]_\gamma$ is the standard $L^\gamma(B_R^{\,\rm c})$ norm, $M$ is a constant depending on $\sup_n\|u_n\|$ but not on $\phi$. A similar inequality for $\int_{B^{\,\rm c}_R}\vert f(x,u_n)\vert\vert\phi\vert$ is also true, therefore
\begin{align}
\int_{B^{\,\rm c}_R}\vert f(x,u_n)-f(x,u)\vert\vert\phi\vert\le\int_{B^{\,\rm c}_R}\vert f(x,u_n)\vert\vert\phi\vert+\int_{B^{\,\rm c}_R}\vert f(x,u)\vert\vert\phi\vert\le 2M\varepsilon\text{.}
\label{xx}
\end{align}
By the compactness of the embedding $X\hookrightarrow L^2_{\rm loc}(\mathbb{R}^2)$ we have
    \begin{equation}
	\sup_{\|\phi\|=1}\int_{B_R}\vert f(x,u_n)-f(x,u)\vert\vert\phi\vert\rightarrow0.\label{eq:42}
    \end{equation}
 It follows from \eqref{xx} and \eqref{eq:42} that
 \begin{align*}
&\hspace{-1cm}\varlimsup_{n\rightarrow\infty}\left\Vert K^{\prime}(u_{n})-K^{\prime
}(u)\right\Vert =\varlimsup_{n\rightarrow\infty}\sup_{\left\Vert
\phi\right\Vert =1}\left\vert \int\left(  f(x,u_{n})-f(x,u)\right)
\phi\right\vert \\
& \leq\varlimsup_{n\rightarrow\infty}\sup_{\left\Vert \phi\right\Vert
=1} \int_{B_{R}}\left|  f(x,u_{n})-f(x,u)\right| \phi
+\varlimsup_{n\rightarrow\infty}\sup_{\left\Vert \phi\right\Vert =1}
\int_{B_{R}^{\,\mathrm{c}}}\left(  f(x,u_{n})-f(x,u)\right)  \phi\\
&\le 2M\varepsilon\text{.}
\end{align*}
Let $\varepsilon\to0$, we deduce $K^{\prime}(u_n)\rightarrow K^{\prime}(u)$ in $X$. Hence $K'$ is compact.
\end{proof}

\begin{rem}\label{rr0}
Lemma \ref{lf:1} is motivated by \cite[Lemma 1]{MR1457116}, In that paper, the space dimension $N>2$ and the working space is $\mathcal{D}^{1,2}(\mathbb{R}^N)$, $a$ and $b$ have to be taken from $L^\infty(\mathbb{R}^N)\cap L^\gamma(\mathbb{R}^N)$ for certain $\gamma$ depending on the power of $|t|$ on the left hand side of \eqref{eq:1.5}. Here our space dimension is $N=2$ and our working space is $X=H^1(\mathbb{R}^2)$. Unlike $\mathcal{D}^{1,2}(\mathbb{R}^N)$ the Sobolev space $X=H^1(\mathbb{R}^2)$ can be continuously embedded into $L^\gamma(\mathbb{R}^2)$ for any $\gamma\ge2$. For this reason, in our assumption $(f_5)$ the integrability of the weight functions $a$ and $b$ can be quit flexible.
\end{rem}

\begin{proof}
	[\indent Proof of Theorem \ref{th:2}]
	As we have pointed out in Remark \ref{re:1}, our condition $(f_5)$ is to ensure the global compactness of our functional $\Phi$. According to the proof of Lemma \ref{l:2.6}, it suffices to derive
	\[
		\varlimsup\int f(x,u_n)(u_n-u)=0.
	\]
	from $u_n\rightharpoonup u$ in $X$.
	
	 Indeed, if $u_n\rightharpoonup u$  in $X$, by the Lemma \ref{lf:1} we have $K^{\prime}(u_n)\rightarrow K^{\prime}(u)$ and
	\begin{align}
		\left\vert \int f(x,u_n)(u_n-u)\right\vert
		&=\left\vert \left\langle K^{\prime}(u_n),u_n-u\right\rangle \right\vert\nonumber\\
		&\leq \left\vert \left\langle K^{\prime}(u_n)-K^{\prime}(u),u_n-u\right\rangle \right\vert+\left\vert \left\langle K^{\prime}(u),u_n-u\right\rangle \right\vert\nonumber\\
		&\leq\left\Vert K^{\prime}(u_n)-K^{\prime}(u)\right\Vert \left\Vert u_n-u\right\Vert+o(1)\rightarrow0.\nonumber
	\end{align}
Therefore, similar to the proof of Lemma \ref{l:2.6} we can show that the functional $\Phi$ satisfies the $(PS)$ condition. Furthermore, applying Lemma \ref{l:3.3} and \ref{l:3.4}, using the same proof we deduce that under the assumptions of Theorem \ref{th:2}, $\Phi$ has a nonzero critical point. This completes the proof of Theorem \ref{th:2}.
\end{proof}
\subsection*{Acknowledgments}
This work was supported by NSFC (12071387)  and NSFC (11971436). The authors would like to thank the anonymous referees for careful reading and valuable suggestions which improve the work.

\end{document}